\documentclass[conference,10pt]{IEEEtran}
\usepackage[dvips]{graphicx}
\usepackage{float}
\usepackage{times}
\usepackage{amsmath}
\usepackage{amsfonts, amsthm}
\renewcommand{\baselinestretch}{1.07}
\usepackage{cite}

\begin{document}

\title{Three-Dimensional Wind Profile Prediction with Trinion-Valued Adaptive Algorithms}
\author{\IEEEauthorblockN{Xiaoming~Gou\IEEEauthorrefmark{1}\IEEEauthorrefmark{2},
Zhiwen~Liu\IEEEauthorrefmark{1},
Wei~Liu\IEEEauthorrefmark{2} and
Yougen~Xu\IEEEauthorrefmark{1}}
\IEEEauthorblockA{\IEEEauthorrefmark{1}School of Information and Electronics, Beijing Institute of Technology, Beijing 100081, China}
\IEEEauthorblockA{\IEEEauthorrefmark{2}Communications Research Group, Department of Electronic and Electrical Engineering,\\ University of Sheffield, Sheffield S1 3JD, United Kingdom\\Email: w.liu@sheffield.ac.uk}}
\maketitle
\begin{abstract}
The problem of three-dimensional (3-D) wind profile prediction is addressed based a trinion wind model, which inherently reckons the coupling of the three perpendicular components of a wind field. The augmented trinion statistics are developed and employed to enhance the prediction performance due to its full exploitation of the second-order statistics. The proposed trinion domain processing can be regarded as a more compact version of the existing quaternion-valued approach, with a lower computational complexity. Simulations based on recorded wind data are provided to demonstrate the effectiveness of the proposed methods.
\end{abstract}

\IEEEpeerreviewmaketitle

\section{Introduction}

Hypercomlex numbers are high-dimensional extensions of real numbers \cite{KantorIL1989,WardJP1997}, and they have been introduced to solve multivariate signal processing problems, such as colour image processing \cite{LeBihanN2003,PeiSC2004,EllTA2007,GuoLQ2011}, vector-sensor array signal processing \cite{MironS2006,LeBihanN2007,GongXF2011,JiangJF2010,ZhangXR2014,liu15a}, human gesture spotting \cite{CheUjangB2014,TobarFA2014}, wind profile prediction \cite{TookCC2011,NavarroMorenoJ2013,BarthelemyQ2014,JiangMD2014}, and wireless communications \cite{liu14n}. In the last application, anemometer readings are modeled with hypercomplex numbers and the wind profile is predicted by adaptive filtering algorithms. In particular, pure quaternions have been widely used to model three-dimensional (3D) wind speed. When external atmospheric parameters are available, a full quaternion-valued model can be considered \cite{TookCC2011}.

The quaternion-valued model leads to improved performance over real-valued models, since it accounts for the coupling of the wind measurements and can be extended to exploiting  the augmented quaternion statistics \cite{TookCC2011}. However,  pure quaternions do not form a mathematical ring \cite{AllenbyRB1991}, since the product of two pure quaternions is not a pure quaternion in general. As a result, the related adaptive algorithms for 3-D wind profile prediction  initialised in the pure quaternion will have to work in the full quaternion domain. The prediction results have to be truncated from full quaternions to pure quaternions, implying redundant computations in the update process. For example, it takes 16 real-valued multiplications and 12 real-valued additions to implement a full-quaternion-valued multiplication. By comparison, it only takes 9 real-valued multiplications and 6 real-valued additions to calculate the multiplication of any two 3-D numbers.

In this work, we aim to develop a more compact 3-D wind speed model based on a 3-D mathematical ring called trinions, termed by an anonymous author who provided all possible definitions of 3-D numbers \cite{Async}. We first define the gradient operation in the trinion domain and then derive the least mean squares (LMS) adaptive algorithm, which is further extended to augmented statistics.

The rest of this paper is organised as follows. A brief introduction to trinions
is provided in Section II. The trinion-valued LMS (TLMS) algorithm is developed in Section III, and  the augmented TLMS (ATLMS) algorithm in Section IV. Simulation results are presented in Section V and conclusions are drawn in Section VI.

\section{Trinions}
As shown in Fig. 1, the 3-D wind speed is a tri-variate signal composed of three perpendicular components, and can be modeled with 3-D hypercomplex numbers (see Fig. 1). There are various definitions of a 3-D number $v$ composed of one real part ($v_a$) and two imaginary parts ($v_b$ and $v_c$),
\begin{equation}
v=v_a+\imath v_b+\jmath v_c,
\end{equation}
according to definitions about the relationships among the three base elements 1, $\imath$, and $\jmath$.
\begin{figure}[t]
  \centering
  \includegraphics[width=1.7in]{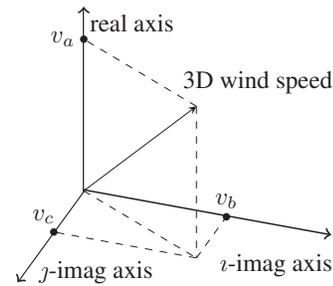}\\
  \caption{Trinion-valued three-dimensional wind speed model.}
\end{figure}

To form an Abelian group of these three, the following rules apply \cite{AssefaD2011}:
\begin{equation}
\imath^2=\jmath,\imath\jmath=\jmath\imath=-1,\jmath^2=-\imath\;.
\end{equation}
Trinions subject to the above rules form a commutative ring, namely, for two trinions $v$ and $w$,  we have  $vw=wv$.

The modulus of $v$ is given by \cite{AssefaD2011}
\begin{equation}
|v|=\sqrt{v_a^2+v_b^2+v_c^2}.
\end{equation}
We define the following conjugate of $v$,
\begin{equation}
v^\ast=v_a-\jmath v_b-\imath v_c,
\end{equation}
so that $|v|^2$ is equal to the real part of $vv^\ast$, i.e. $|v|^2=\Re(vv^\ast)$.

To develop an LMS-like algorithm in the trinion domain, the trinion-valued gradient needs to be defined. In the complex domain, a variable $z$ and its conjugate can be viewed as two independent variables, based on which the complex-valued gradient can be defined \cite{vandenBosA1994}. Similarly, in the quaternion domain, a variable $q$ and its three involutions can be viewed as four independent variables, based on which the quaternion-valued gradient can be derived \cite{JiangMD2014,liu14p}. However, to our best knowledge, the trinion involution is not available in general. In this paper, we define the following gradients of a function $f(\boldsymbol x)$ with respect to the trinion-valued variable $\boldsymbol x$ and its conjugate,
\begin{equation}
\begin{split}
\nabla_{\boldsymbol x}f&=\frac{1}{3}\big(\nabla_{\boldsymbol x_a}f-\jmath\nabla_{\boldsymbol x_b}f-\imath\nabla_{\boldsymbol x_c}f\big),\\
\nabla_{\boldsymbol x^\ast}f&=\frac{1}{3}\big(\nabla_{\boldsymbol x_a}f+\imath\nabla_{\boldsymbol x_b}f+\jmath\nabla_{\boldsymbol x_c}f\big),
\end{split}
\end{equation}
respectively, where $\boldsymbol x=\boldsymbol x_a+\imath\boldsymbol x_b+\jmath\boldsymbol x_c$.

Unlike the quaternion-valued gradient, the trinion-valued gradient remains the same, no matter which side of the sub-gradients the imaginary units $\imath$ and $\jmath$ are on, since trinions are commutative. We can then calculate the derivatives of some simple functions, for instance,
\begin{equation}
\begin{split}
\frac{\partial x}{\partial x}=\frac{\partial x^\ast}{\partial x^\ast}&=1,\\
\frac{\partial x}{\partial x^\ast}=\frac{\partial x^\ast}{\partial x}&=\frac{1-\imath+\jmath}{3}.
\end{split}
\end{equation}

\section{Trinion-valued LMS Filtering}
Similar to the LMS algorithm in the complex-valued domain, the trinion-valued error $e(n)$ is given by
\begin{equation}
e(n)=d(n)-\boldsymbol w^\mathrm T(n)\boldsymbol x(n),
\end{equation}
where $d(n)$ is the reference signal, $\boldsymbol w(n)$ is the weight vector, $\boldsymbol x(n)=[x(n-P);x(n-P-1);\cdots;x(n-L-P +1)]^T$ is the filter input,
$L$ is the filter length, and $P$ is the prediction step. The cost function is expressed as
\begin{equation}
J(n)=|e(n)|^2.
\end{equation}
Using the steepest descent method, the following gradient need to be calculated
\begin{equation}
\nabla_{\boldsymbol w^\ast}J(n)=\frac{1}{3}\big[\nabla_{\boldsymbol w_a}J(n)+\imath\nabla_{\boldsymbol w_b}J(n)+\jmath\nabla_{\boldsymbol w_c}J(n)\big].\label{gra_def}
\end{equation}
For details of the calculation, please refer to Appendix at the end. Using results there, the update of the weight vector can be obtained as
\begin{equation}
\boldsymbol w(n+1)=\boldsymbol w(n)+\mu e(n)\boldsymbol x^\ast(n),\label{update}
\end{equation}
where $\mu$ is the step size. The algorithm formulated above is termed as the trinion-valued LMS (TLMS) algorithm.

\section{Augmented Trinion Statistics}
A zero-mean trinion-valued vector $\boldsymbol v$ is composed of three zero-mean real-valued random vectors $\boldsymbol v_a$, $\boldsymbol v_b$, and $\boldsymbol v_c$, and their complete second-order statistics can be found in the following six real-valued covariance matrices:
\begin{equation}
\begin{split}
\boldsymbol C_{\boldsymbol v_a\boldsymbol v_a}=E\{\boldsymbol v_a\boldsymbol v_a^\mathrm T\},\quad\boldsymbol C_{\boldsymbol v_b\boldsymbol v_b}=E\{\boldsymbol v_b\boldsymbol v_b^\mathrm T\},\\
\quad\boldsymbol C_{\boldsymbol v_c\boldsymbol v_c}=E\{\boldsymbol v_c\boldsymbol v_c^\mathrm T\},\quad\boldsymbol C_{\boldsymbol v_a\boldsymbol v_b}=E\{\boldsymbol v_a\boldsymbol v_b^\mathrm T\},\\
\quad\boldsymbol C_{\boldsymbol v_b\boldsymbol v_c}=E\{\boldsymbol v_b\boldsymbol v_c^\mathrm T\},\quad\boldsymbol C_{\boldsymbol v_c\boldsymbol v_a}=E\{\boldsymbol v_c\boldsymbol v_a^\mathrm T\}.
\end{split}
\end{equation}
And they can be more efficiently represented by three trinion-valued covariance matrices:
\begin{equation}
\begin{split}
\boldsymbol C_{\boldsymbol v\boldsymbol v}&=E\{\boldsymbol v\boldsymbol v^\mathrm H\},\\
\boldsymbol C_{\boldsymbol v\boldsymbol v^\imath}&=E\{\boldsymbol v\boldsymbol v^{\imath\mathrm H}\},\\
\boldsymbol C_{\boldsymbol v\boldsymbol v^\jmath}&=E\{\boldsymbol v\boldsymbol v^{\jmath\mathrm H}\},
\end{split}
\end{equation}
where superscript $^\mathrm H$ denotes the Hermitian transpose and the two additional mappings of $\boldsymbol v$ are defined as
\begin{equation}
\begin{split}
\boldsymbol v^\imath&=\boldsymbol v_b-\imath\boldsymbol v_a-\jmath\boldsymbol v_c,\\
\boldsymbol v^\jmath&=\boldsymbol v_c-\imath\boldsymbol v_b-\jmath\boldsymbol v_a.
\end{split}
\end{equation}
Neither of these two mappings is an involution, and they are defined as shorthand notations only. Then we can obtain
\begin{equation}
\begin{split}
\boldsymbol C_{\boldsymbol v_a\boldsymbol v_a}&=\frac{1}{2}\Re(\boldsymbol C_{\boldsymbol v\boldsymbol v}+\jmath\boldsymbol C_{\boldsymbol v\boldsymbol v^\imath}),\\
\boldsymbol C_{\boldsymbol v_b\boldsymbol v_b}&=\frac{1}{2}\Re(\jmath\boldsymbol C_{\boldsymbol v\boldsymbol v^\imath}-\jmath\boldsymbol C_{\boldsymbol v\boldsymbol v^\imath}),\\
\boldsymbol C_{\boldsymbol v_c\boldsymbol v_c}&=\frac{1}{2}\Re(\boldsymbol C_{\boldsymbol v\boldsymbol v}-\imath\boldsymbol C_{\boldsymbol v\boldsymbol v^\jmath}),\\
\boldsymbol C_{\boldsymbol v_a\boldsymbol v_b}&=\frac{1}{2}\Re(\boldsymbol C_{\boldsymbol v\boldsymbol v^\imath}+\jmath\boldsymbol C_{\boldsymbol v\boldsymbol v^\jmath}),\\
\boldsymbol C_{\boldsymbol v_b\boldsymbol v_c}&=\frac{1}{2}\Re(\imath\boldsymbol C_{\boldsymbol v\boldsymbol v}-\jmath\boldsymbol C_{\boldsymbol v\boldsymbol v^\jmath}),\\
\boldsymbol C_{\boldsymbol v_c\boldsymbol v_a}&=\frac{1}{2}\Re(\boldsymbol C_{\boldsymbol v\boldsymbol v^\imath}-\imath\boldsymbol C_{\boldsymbol v\boldsymbol v}).
\end{split}
\end{equation}
To take the complete second-order statistics into consideration, we need to make use of the augmented input vector composed of the original input vector $\boldsymbol x(n)$ and its two mappings $\boldsymbol x^\imath(n),\boldsymbol x^\jmath(n)$, namely,
\begin{equation}
\boldsymbol x^{\mathrm{aug}}(n)=
\begin{bmatrix}
\boldsymbol x(n)\\
\boldsymbol x^\imath(n)\\
\boldsymbol x^\jmath(n)
\end{bmatrix},
\end{equation}
and the predicted estimate $y(n)$ will be given by
\begin{equation}
\begin{split}
y(n)&=\boldsymbol w^\mathrm{augT}(n)\boldsymbol x^{\mathrm{aug}}(n)\\
&=\boldsymbol w_1^\mathrm T(n)\boldsymbol x(n)+\boldsymbol w_2^\mathrm T(n)\boldsymbol x^\imath(n)+\boldsymbol w_3^\mathrm T(n)\boldsymbol x^\jmath(n),
\end{split}
\end{equation}
where $\boldsymbol w^\mathrm{aug}=[\boldsymbol w_1;\boldsymbol w_2;\boldsymbol w_3]$. Analogous to the derivation of the TLMS algorithm, the update equation of the augmented TLMS (ATLMS) algorithm is given by
\begin{equation}
\boldsymbol w^\mathrm{aug}(n+1)=\boldsymbol w^\mathrm{aug}(n)+\rho e^\mathrm{aug}(n)\boldsymbol x^{\mathrm{aug}\ast}(n),\label{aug_update}
\end{equation}
where $\rho$ is the step size.

The computational complexity for each  update of the weight vector of the trinion-based and quaternion-based filtering algorithms are summarised in Table I. It can be seen that the trinion model can effectively reduce the computational complexity compared to the quaternion model.
\begin{table}[h]
\caption{Computional complexity per update of the weight vector}
\centering
\begin{tabular}{c|c|c}
\hline
Algorithm & Real multiplications & Real additions\\
\hline
QLMS \cite{BarthelemyQ2014}& $16L+4$ & $16L$\\
Augmented QLMS & $64L+4$ & $64L$\\
TLMS & $9L+3$ & $9L$\\
Augmented TLMS & $27L+3$ & $27L$\\
\hline
\end{tabular}
\end{table}

\begin{figure}
  \centering
  \includegraphics[width=3.7in]{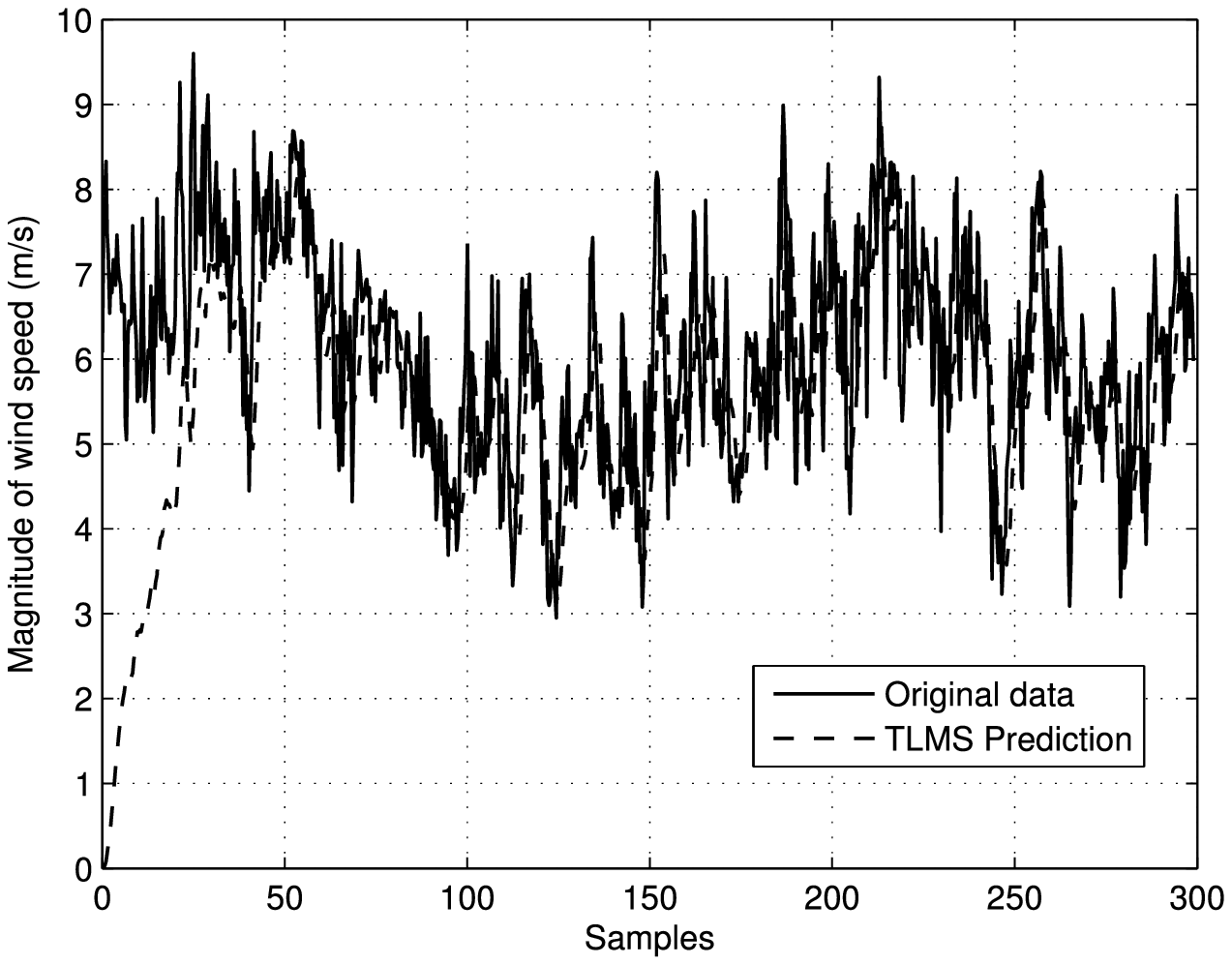}\\
  \caption{Predicted result from the TLMS algorithm.}
\end{figure}

\begin{figure}
  \centering
  \includegraphics[width=3.7in]{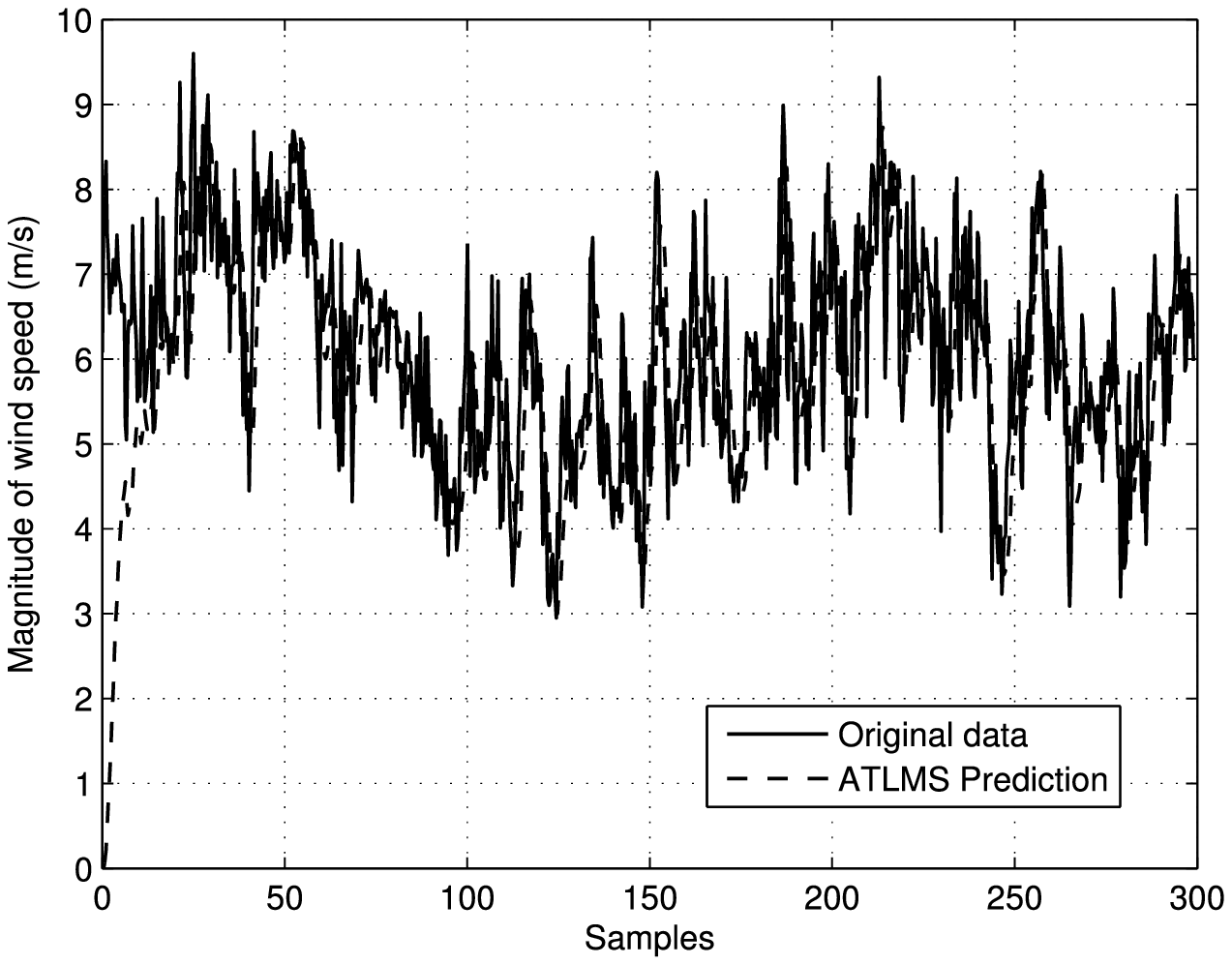}\\
  \caption{Predicted result from the ATLMS algorithm.}
\end{figure}

\begin{figure}
  \centering
  \includegraphics[width=3.7in]{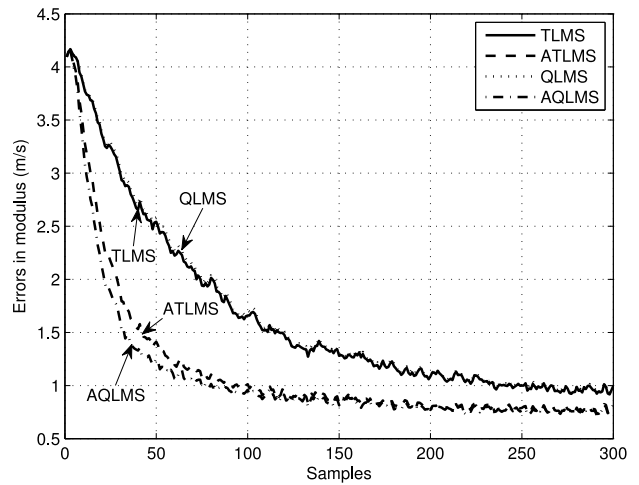}\\
  \caption{Averaged learning curves.}
\end{figure}

\section{Simulations}
In this section, both proposed algorithms (TLMS and ATLMS) are applied to anemometer readings provided by Google's RE$<$C Initiative \cite{Google}. The wind speed measured on May 25, 2011 is used for demonstration. The step size is set to be $6\times10^{-5}$. The filter length is 8, and the prediction step is 1. All algorithms are initialised with zeros. The predicted results provided by TLMS and ATLMS algorithms are shown in Figs. 2 and  3, respectively. From the results, we can see that both algorithms can track the wind data effectively.

The learning curves averaged over 200 trials of the proposed algorithms are shown in Fig. 4, compared with the quaternion-based QLMS and AQLMS algorithms. It can be observed that both augmented algorithms (AQLMS and ATLMS) have a faster convergence rate than the original ones (QLMS and TLMS), due to their full exploitation of the second-order statistics. Meanwhile, the proposed TLMS algorithm has a similar performance with the QLMS algorithm, while the ATLMS algorithm is comparable with the AQLMS algorithm. However, as shown in Table I, the proposed trinion-based algorithms have a much lower computational complexity.

\section{Conclusion}
A compact model for three-dimensional wind profile prediction based on  trinion algebra has been proposed, with two LMS-type adaptive filtering algorithms derived using partial and full trinion-valued second-order statistics, respectively. Numerical simulations using recorded wind data have shown that the two algorithms have a similar  performance to their quaternion-valued counterparts, but   with  a much lower computational complexity.

\section*{Appendix}
The cost function $J(n)$ as a function of real-valued variables is given by
\begin{equation}
\begin{split}
J&=(d_a-\boldsymbol w_a^\mathrm T\boldsymbol x_a+\boldsymbol w_b^\mathrm T\boldsymbol x_c+\boldsymbol w_c^\mathrm T\boldsymbol x_b)^2\\
&+(d_b-\boldsymbol w_a^\mathrm T\boldsymbol x_b-\boldsymbol w_b^\mathrm T\boldsymbol x_a+\boldsymbol w_c^\mathrm T\boldsymbol x_c)^2\\
&+(d_c-\boldsymbol w_a^\mathrm T\boldsymbol x_c-\boldsymbol w_b^\mathrm T\boldsymbol x_b-\boldsymbol w_c^\mathrm T\boldsymbol x_a)^2,
\end{split}
\end{equation}
where the time index `$n$' has been dropped for the sake of compact notation. Then the three component-wise gradients can be computed as
\begin{equation}
\begin{split}
\nabla_{\boldsymbol w_a}J&=2\big[(\boldsymbol x_a\boldsymbol x_a^\mathrm T+\boldsymbol x_b\boldsymbol x_b^\mathrm T+\boldsymbol x_c\boldsymbol x_c^\mathrm T)\boldsymbol w_a\\
&+(\boldsymbol x_b\boldsymbol x_a^\mathrm T+\boldsymbol x_c\boldsymbol x_b^\mathrm T-\boldsymbol x_a\boldsymbol x_c^\mathrm T)\boldsymbol w_b\\
&+(\boldsymbol x_c\boldsymbol x_a^\mathrm T-\boldsymbol x_a\boldsymbol x_b^\mathrm T-\boldsymbol x_b\boldsymbol x_c^\mathrm T)\boldsymbol w_c\\
&-(d_a\boldsymbol x_a+d_b\boldsymbol x_b+d_c\boldsymbol x_c)\big],
\end{split}
\label{sub_gra1}
\end{equation}
\begin{equation}
\begin{split}
\nabla_{\boldsymbol w_b}J&=2\big[(\boldsymbol x_a\boldsymbol x_b^\mathrm T+\boldsymbol x_b\boldsymbol x_c^\mathrm T-\boldsymbol x_c\boldsymbol x_a^\mathrm T)\boldsymbol w_a\\
&+(\boldsymbol x_c\boldsymbol x_c^\mathrm T+\boldsymbol x_a\boldsymbol x_a^\mathrm T+\boldsymbol x_b\boldsymbol x_b^\mathrm T)\boldsymbol w_b\\
&+(\boldsymbol x_c\boldsymbol x_b^\mathrm T-\boldsymbol x_a\boldsymbol x_c^\mathrm T+\boldsymbol x_b\boldsymbol x_a^\mathrm T)\boldsymbol w_c\\
&+(d_a\boldsymbol x_c-d_b\boldsymbol x_a-d_c\boldsymbol x_b)\big],
\end{split}
\label{sub_gra2}
\end{equation}
\begin{equation}
\begin{split}
\nabla_{\boldsymbol w_c}J&=2\big[(\boldsymbol x_a\boldsymbol x_c^\mathrm T-\boldsymbol x_b\boldsymbol x_a^\mathrm T-\boldsymbol x_c\boldsymbol x_b^\mathrm T)\boldsymbol w_a\\
&+(\boldsymbol x_b\boldsymbol x_c^\mathrm T-\boldsymbol x_c\boldsymbol x_a^\mathrm T+\boldsymbol x_a\boldsymbol x_b^\mathrm T)\boldsymbol w_b\\
&+(\boldsymbol x_a\boldsymbol x_a^\mathrm T+\boldsymbol x_c\boldsymbol x_c^\mathrm T+\boldsymbol x_b\boldsymbol x_b^\mathrm T)\boldsymbol w_c\\
&+(d_a\boldsymbol x_b+d_b\boldsymbol x_c-d_c\boldsymbol x_a)\big].
\end{split}
\label{sub_gra3}
\end{equation}
Finally, we can obtain the expression for the gradient of $J(n)$ by substituting (\ref{sub_gra1})--(\ref{sub_gra3}) into (\ref{gra_def}),
\begin{equation}
\nabla_{\boldsymbol w^\ast}J(n)=\frac{2}{3}e(n)\boldsymbol x^\ast(n),
\end{equation}
which yields the update equation in (\ref{update}), as the term $\frac{2}{3}$ can be absorbed into the step size.

\section*{Acknowledgment}
This work was supported in part by the National Natural Science Foundation of China (61331019, 61490690), the China Scholarship Council Postgraduate Scholarship Program (2014), and the UK National Grid.



\begin{thebibliography}{20}

\bibitem{KantorIL1989}
I. L. Kantor and A. S. Solodovnikov, \emph{Hypercomplex numbers: An elementary introduction to algebras}. Springer-Verlag, 1989.
\bibitem{WardJP1997}
J. P. Ward, \emph{Quaternions and Cayley numbers: Algebra and applications}. Dordrecht: Kluwer Academic Publisher, 1997.
\bibitem{LeBihanN2003}
N. Le Bihan and S. J. Sangwine, ``Quaternion principle component analysis of color images,'' in \emph{Proceedings of the IEEE International Conference on Image Processing}, Barcelona, Spain, pp. 809--812, Sep. 2003.
\bibitem{PeiSC2004}
S. C. Pei, J. H. Chang, and J. J. Ding, ``Commutative reduced biquaternions and their Fourier transform for signal and image processing applications,'' \emph{IEEE Transactions on Signal Processing}, vol. 52, no. 7, pp. 2012--2031, Jul. 2004.
\bibitem{EllTA2007}
T. A. Ell, ``Multi-vector color-image filters,'' in \emph{Proceedings of the IEEE International Conference on Image Processing}, San Antonio, TX, USA, pp. 245--248, Sep. 2007.
\bibitem{GuoLQ2011}
L. Q. Guo, M. Zhu, and X. H. Ge, ``Reduced biquaternion canonical transform, convolution and correlation,'' \emph{Signal Processing}, vol. 91, no. 8, pp. 2147--2153, Aug. 2011.
\bibitem{MironS2006}
S. Miron, N. Le Bihan, and J. I. Mars, ``Quaternion-MUSIC for vector-sensor array processing,'' \emph{IEEE Transactions Signal Processing}, vol. 54, no. 4, pp. 1218--1229, Apr. 2006.
\bibitem{LeBihanN2007}
N. Le Bihan, S. Miron, and J. I. Mars, ``MUSIC algorithm for vector-sensors array processing using biquaternions,'' \emph{IEEE Transactions on
Signal Processing}, vol. 55, no. 9, pp. 4523--4533, Sep. 2007.
\bibitem{JiangJF2010}
J. F. Jiang and J. Q. Zhang, ``Geometric algebra of Euclidean 3-space for electromagnetic vector-sensor array processing, part I: Modeling,'' \emph{IEEE Transactions on Antennas and Propagation}, vol. 58, no. 12, pp. 3961--3973, Dec. 2010.
\bibitem{GongXF2011}
X. F. Gong, Z. W. Liu, and Y. G. Xu, ``Coherent source localization: Bicomplex polarimetric smoothing with electromagnetic vector-sensors,'' \emph{IEEE Transactions on Aerospace and Electronic Systems}, vol. 47, no. 3, pp. 2268--2285, Jul. 2011.
\bibitem{ZhangXR2014}
X. R. Zhang, W. Liu, Y. G. Xu, and Z. W. Liu, ``Quaternion-valued robust adaptive beamformer for electromagnetic vector-sensor arrays with worst-case constraint,'' \emph{Signal Processing}, vol. 104, pp. 274--283, Nov. 2014.
\bibitem{liu15a}
M.~B. Hawes and W.~Liu, \newblock ``Design of fixed beamformers based on vector-sensor arrays,''
\newblock {\em International Journal of Antennas and Propagation}, vol. 2015, 2015.
\bibitem{CheUjangB2014}
B. Che Ujang, C. Jahanchahi, C. C. Took, and D. P. Mandic, ``Adaptive convex combination approach for the identification of improper quaternion processes,'' \emph{IEEE Transactions on Neural Networks and Learning Systems}, vol. 25, no. 1, pp. 172--182, Jan. 2014.
\bibitem{TobarFA2014}
F. A. Tobar and D. P. Mandic, ``Quaternion reproducing kernel Hilbert spaces: Existence and uniqueness conditions,'' \emph{IEEE Transactions on Information Theory}, vol. 60, no. 9, pp. 5736--5749, Sep. 2014.
\bibitem{TookCC2011}
C. C. Took, G. Strbac, K. Aihara, and D. P. Mandic, ``Quaternion-valued short term joint forecasting of three-dimensional wind and atmospheric parameters,'' \emph{Renewable Energy}, vol. 36, no. 6, pp. 1754--1760, Jun. 2011.
\bibitem{NavarroMorenoJ2013}
J. Navarro-Moreno, R. M. Fernandez-Alcal\'a, C. C. Took, and D. P. Mandic, ``Prediction of wide-sense stationary quaternion random signals,'' \emph{Signal Processing}, vol. 93, no. 9, pp. 2573--2580, Sep. 2013.
\bibitem{BarthelemyQ2014}
Q. Barth{\'e}lemy, A. Larue, and J. I. Mars, ``About QLMS derivations,'' \emph{IEEE Signal Processing Letters}, vol. 21, no. 2, pp. 240--243, Feb. 2014.
\bibitem{AllenbyRB1991}
R. B. Allenby. \emph{Rings, fields and groups: An introduction to abstract algebra}. Butterworth-Heinemann, 1991.
\bibitem{Async}
Async, \emph{OMIC's N-nion's site} [Online]. Available: http://asyncbrain.baf.cz/m/nt/index.htm.
\bibitem{AssefaD2011}
D. Assefa, L. Mansinha, K. F. Tiampo, H. Rasmussen, and K. Abdella, ``The trinion Fourier transform of color images,'' \emph{Signal Processing}, vol. 91, no. 8, pp. 1887--1900, Aug. 2011.
\bibitem{vandenBosA1994}
A. van den Bos, ``Complex gradient and Hessian,'' \emph{IEE Proceedings - Vision Image Signal Process}, vol. 141, no. 6, pp. 380--383, Dec. 1994.
\bibitem{JiangMD2014}
M. D. Jiang, W. Liu, and Y. Li, ``A general quaternion-valued gradient operator and its applications to computational fluid dynamics and adaptive beamforming,'' in \emph{Proceedings of the 19th International Conference on Digital Signal Processing}, Hong Kong, China, pp. 821--826, Aug. 2014.
\bibitem{liu14n}
W.~Liu,
\newblock ``Antenna array signal processing for a quaternion-valued wireless
  communication system,''
\newblock in {\em Proc. the Benjamin Franklin Symposium on Microwave and
  Antenna Sub-systems (BenMAS)}, Philadelphia, US, September 2014.
\bibitem{liu14p}
M.~D. Jiang, Y.~Li, and W.~Liu,
\newblock ``Properties and applications of a restricted \mbox{HR} gradient
  operator,''
\newblock {\em arXiv:1407.5178 [math.OC]}, July 2014.
\bibitem{Google}
Google, ``RE$<$C: Surface level wind data collection,'' \emph{Google Code 2011} [Online]. Available: http://code.google.com/p/google-rec-csp/
\end{thebibliography}
\end{document}